\documentclass{amsart}
\usepackage{amssymb} 
\usepackage{enumitem}
\usepackage{multirow}
\setlist[itemize]{leftmargin=12mm}
\setlist[enumerate]{leftmargin=12mm}
\usepackage{comment}

\DeclareMathOperator{\Aut}{Aut}
\DeclareMathOperator{\Cl}{Cl}

\DeclareMathOperator{\GL}{GL}

\DeclareMathOperator{\SL}{SL}

\DeclareMathOperator{\Gal}{Gal}
\DeclareMathOperator{\norm}{Norm}
\DeclareMathOperator{\ord}{\upsilon}
\newcommand{\Q}{{\mathbb Q}}
\newcommand{\Z}{{\mathbb Z}}

\newcommand{\F}{{\mathbb F}}

\newcommand{\cA}{\mathcal{A}}

\newcommand{\cG}{\mathcal{G}}

\newcommand{\cH}{\mathcal{H}}

\newcommand{\cM}{\mathcal{M}}
\newcommand{\cN}{\mathcal{N}}

\newcommand{\OO}{{\mathcal O}}
\newcommand{\ga}{\mathfrak{a}}
\newcommand{\gb}{{\mathfrak{b}}}

\newcommand{\ff}{\mathfrak{f}}

\newcommand{\gp}{\mathfrak{p}}
\newcommand{\fp}{\mathfrak{m}}
\newcommand{\fr}{\mathfrak{r}}
\newcommand{\fq}{\mathfrak{q}}
\newcommand{\mP}{\mathfrak{P}}

\begin {document}

\newtheorem{thm}{Theorem}
\newtheorem{lem}{Lemma}[section]
\newtheorem{prop}[lem]{Proposition}

\newtheorem{cor}[lem]{Corollary}

\theoremstyle{definition}

\theoremstyle{remark}

\title[Fermat]{
Fermat's Last Theorem over some Small Real Quadratic Fields
}
\author{Nuno Freitas and Samir Siksek}
\address{
 Mathematisches Institut\\
Universit\"{a}t Bayreuth\\
95440 Bayreuth, Germany
}

\address{Mathematics Institute\\
	University of Warwick\\
Coventry\\
	CV4 7AL \\
	United Kingdom}


\date{\today}
\thanks{
The first-named author is supported
through a grant within the framework of the DFG Priority Programme 1489
{\em Algorithmic and Experimental Methods in Algebra, Geometry and Number Theory}.
The second-named
author is supported by an EPSRC Leadership Fellowship EP/G007268/1,
and EPSRC {\em LMF: L-Functions and Modular Forms} Programme Grant
EP/K034383/1.
}

\keywords{Fermat, modularity, Galois representation, level lowering}
\subjclass[2010]{Primary 11D41, Secondary 11F80, 11F03}

\begin{abstract}
Using modularity, level lowering, and explicit
computations with Hilbert modular forms, Galois representations and ray class groups, we show
that for $3 \le d \le 23$ squarefree, $d \ne 5$, $17$,
the Fermat equation $x^n+y^n=z^n$ has no non-trivial solutions
over the quadratic field $\Q(\sqrt{d})$ for $n \ge 4$. 
Furthermore, we show for $d=17$ that the same holds for prime
exponents $n \equiv 3$, $5 \pmod{8}$.
\end{abstract}


\maketitle



\section{Introduction}

Interest in the Fermat equation
\begin{equation}\label{eqn:Fermatn}
x^n+y^n=z^n, 
\end{equation}
 over various number
fields goes back to the 19th and early 20th century,
with the work of 
Maillet (1897) and Furtw\"{a}ngler (1910)
(mentioned in Dickson's {\em History of the Theory of Numbers}
\cite[pages 758 and 768]{Dickson}). However, until the current
work, the only number fields for which Fermat's Last Theorem
has been proved are $\Q$ and $\Q(\sqrt{2})$; the former in
the ground breaking work of  
Wiles \cite{Wiles}, and the latter by Jarvis and Meekin \cite{JMee}.
The precise statements are that if $K=\Q$ and $n \ge 3$
or $K=\Q(\sqrt{2})$ and $n \ge 4$, then the only solutions to
\eqref{eqn:Fermatn} in $K$ are the trivial ones satisfying $xyz=0$.
In \cite{Fermat} it is shown for five-sixths of real quadratic fields
$K$ that there is a bound $B_K$, such that for prime exponents $n \ge B_K$ the only
solutions to the Fermat equation \eqref{eqn:Fermatn} over $K$ are
the trivial ones.
This paper is concerned with proving Fermat's Last Theorem for
several other real quadratic fields. Our main results are the following two theorems.
\begin{thm}\label{thm:main}
Let $3 \le d \le 23$ squarefree, $d \ne 5$, $17$. Then the Fermat
equation~\eqref{eqn:Fermatn}
 does not have any non-trivial solutions over $\Q(\sqrt{d})$
with exponent $n \ge 4$.
\end{thm}
\begin{thm}\label{thm:d=17}
The Fermat equation~\eqref{eqn:Fermatn} 
has no non-trivial solutions 
over $\Q(\sqrt{17})$
for prime exponents $n \geq 5$ satisfying
$n \equiv 3$, $5 \pmod{8}$.
\end{thm}
\noindent \textbf{Remark.}
For $n=3$ equation~\eqref{eqn:Fermatn} defines an elliptic curve having rank $0$ over $\Q$;
it does however have positive rank over some of
the quadratic fields in the statement of
Theorem~\ref{thm:main}. We therefore impose $n \geq 4$.

\bigskip

It is sufficient to prove Theorem~\ref{thm:main} for exponents $n=4$, $6$, $9$,
and for prime exponents $n=p \ge 5$. In fact,
 all solutions to the Fermat equation
in quadratic fields for $n=4$, $6$, $9$ have been determined by
Aigner \cite{Aigner1}, \cite{Aigner2}. These are all defined over imaginary
quadratic fields except for the trivial solutions. We may therefore
restrict our attention to prime exponents $n=p \ge 5$.

Let $d \ge 2$ be a squarefree positive integer, and let $K=\Q(\sqrt{d})$,
and write $\OO_K$
for its ring of integers.
By the \textbf{Fermat equation with exponent $p$ over $K$} we mean the equation
\begin{equation}\label{eqn:Fermat}
a^p+b^p+c^p=0, \qquad a,b,c\in \OO_K.
\end{equation}
A solution $(a,b,c)$ is called \textbf{trivial}
if $abc=0$, otherwise \textbf{non-trivial}.
For
$p=5$, $7$, $11$, all solutions of degree $\le (p-1)/2$ have
been determined by Gross and Rohrlich \cite[Theorem 5]{GR}. It turns out that the only
non-trivial quadratic solutions are permutations of $(1,\omega,\omega^2)$,
where $\omega$ is a primitive cube root of unity.
For $p=13$, the same was shown to be true by Tzermias \cite{Tzermias}. We shall
therefore henceforth
assume that $p\ge 17$.

\subsection{A brief overview of the method and difficulties}
As in the proof of Fermat's Last Theorem over $\Q$ and $\Q(\sqrt{2})$,
let $(a,b,c)$ be a non-trivial solution to the Fermat equation~\eqref{eqn:Fermat},
and consider the {\text Frey} elliptic curve
\begin{equation}\label{eqn:Frey}
E_{a,b,c} \; : \; Y^2=X(X-a^p)(X+b^p) .
\end{equation}
Write $E=E_{a,b,c}$ and denote by $\overline{\rho}_{E,p}$  
its mod $p$ Galois representation. 
An essential fact to the proof of Fermat's Last Theorem over $\Q$ and
$\Q(\sqrt{2})$ is the modularity of the Frey curve. Modularity of elliptic
curves over all real quadratic fields is now known (see \cite{FHS}). 
In
particular, our Frey curve $E_{a,b,c}$ is modular over $K$. The proof of Fermat's Last Theorem
over $\Q$ and $\Q(\sqrt{2})$ makes essential use of the fact that it is always possible
to scale and permute the hypothetical non-trivial solution so that not only are $a$, $b$, $c$ 
algebraic integers, but they are also coprime, and they satisfy additional $2$-adic restrictions;
over $\Q$ these are $a \equiv -1 \pmod{4}$ and $2 \mid b$. 
In both cases, a suitable
choice of scaling produces a semistable Frey curve $E_{a,b,c}$. Applying suitable
level-lowering results to the modular Galois representation $\overline{\rho}_{E,p}$
shows that it arises from an eigenform of level $2$ for $\Q$, and a Hilbert
eigenform of level $\sqrt{2}$ for $\Q(\sqrt{2})$. There are no eigenforms
at these levels, giving a contradiction and completing the proof for both fields.

\medskip

It should be possible to carry out the same level-lowering strategy
over any real quadratic field $K$.
To build on  this and prove Fermat's Last Theorem over $K$ there are however three principal difficulties.
\begin{enumerate}
\item[(a)] Verifying the irreducibility of $\overline{\rho}_{E,p}$;
this is required for applying level-lowering theorems.
\item[(b)] Computing the newforms at the levels predicted
by conductor computations and level-lowering; in general
these levels depend on the class and unit groups of $K$ and
are not of small norm.
\item[(c)] Dealing with the newforms found at these levels;
in general these spaces will be non-zero.  
\end{enumerate}
In \cite{Fermat}, it is shown that these difficulties disappear
for $p>C_K$, where $C_K$ is some inexplicit constant, for five-sixths
of real quadratic fields $K$. In this paper we show how to 
deal with (a), (b), (c) for the fields in the statement of Theorem~\ref{thm:main}.
For $K=\Q(\sqrt{d})$ with $d=5$, $d=17$, our method for (c) fails.
However, for $d=17$, we are still able to prove Fermat's Last Theorem
for $1/2$ of exponents $p$ using an arguments of Halberstadt and Kraus \cite{HK},
yielding Theorem~\ref{thm:d=17}.

\medskip

For $K=\Q$, it follows from Mazur's celebrated theorem on isogenies \cite{Mazur} that
$\overline{\rho}_{E,p}$ is irreducible for the Frey curve $E$ 
and $p \ge 5$.  The analogue of Mazur's Theorem is not known
for any other number field. However, Kraus shows
that for $K$ a real quadratic field of class number $1$, and $E$
a semistable elliptic curve over $K$, if $\overline{\rho}_{E,p}$
is reducible, then $p \le 13$ or ramifies in $K$ or $p \mid \norm_{K/\Q}(u^2-1)$
where $u$ is a fundamental unit of $K$. It is possible to give
corresponding  bounds for $p$ 
when the  class group is non-trivial, and  
$E$ has some fixed set of additive primes (see
for example \cite{DavidI}), although these bounds are considerably worse.
In Section~\ref{sec:irred} we overcome these difficulties for
the fields in Theorem~\ref{thm:main} through explicit class field
theory computations. Without this, we would not have been able
to deal with small exponents $p$.

\medskip

Assuming modularity of the Frey curve, it should be possible 
to apply the same strategy to the Fermat equation over 
many totally real fields, although the computation of newforms
for totally real fields with either large degree or large discriminant is likely to be impractical. In Section~\ref{sec:end} we illustrate this
by looking at $\Q(\sqrt{30})$ and $\Q(\sqrt{79})$. We also
include at the end of that section a comparison between
the recipes of the current paper, and those of \cite{Fermat},
illustrating the need for the improvements in the current work.

\medskip

The computations in this paper were carried out on the computer
algebra system {\tt Magma} \cite{magma}. In particular, we have 
used {\tt Magma}'s Hilbert Modular Forms package  
(for the theory see the work of Demb\'{e}l\'{e}, Donnelly and Voight \cite{DD}, \cite{DV})
and Class Field Theory package (due to C.\ Fieker).

\subsection{Notational conventions}\label{sec:conv}
Throughout $p$ denotes an odd rational prime, and $K$ a 
totally real number field, with 
ring of integers $\OO_K$. We write $S$ for the set
of prime ideals in $\OO_K$ dividing 2.
For a non-zero ideal $\ga$ of $\OO_K$, we denote by
$[\ga]$ the class of $\ga$ in the class group $\Cl(K)$. 
For a non-trivial solution $(a,b,c)$ to the Fermat equation \eqref{eqn:Fermat},
let
\begin{equation}\label{eqn:cG}
\cG_{a,b,c}:=a \OO_K+b\OO_K+c\OO_K \, .
\end{equation}
and let $[a,b,c]$ denote
the class of $\cG_{a,b,c}$ in $\Cl(K)$.
We exploit the well-known fact
(e.g. \cite[Theorem VIII.4]{CassFr})
that every ideal class contains
infinitely many prime ideals. 
Let $r=\# (\Cl(K)/\Cl(K)^2)$. 
Let $\fp_1=1 \cdot \OO_K$, and 
fix once and for all odd prime ideals $\fp_2,\dots,\fp_r$ such 
that $[\fp_1],\dots,[\fp_r]$
represent the cosets of $\Cl(K)/\Cl(K)^2$. Let
\[
\cH=\{ \fp_1,\dots,\fp_r\}.
\]

Let $G_K=\Gal(\overline{K}/K)$.
For an elliptic curve $E/K$,
we write
\begin{equation}\label{eqn:rho}
\overline{\rho}_{E,p}\; :\; G_K \rightarrow \Aut(E[p]) \cong \GL_2(\F_p),
\end{equation}
for the representation of $G_K$ on the $p$-torsion of $E$. 
For
a Hilbert eigenform $\ff$ over $K$, we let $\Q_\ff$ denote
the field generated by its eigenvalues.
In this situation $\varpi$ will denote a prime of $\Q_\ff$
above $p$; of course if $\Q_\ff=\Q$ we write
$p$ instead of $\varpi$. All other primes we consider are primes
of $K$. 
We reserve the symbol 
$\mP$ for  primes belonging to $S$, and 
$\fp$ for 
primes belonging to $\cH$.
An arbitrary prime of $K$ is denoted by $\fq$,  and $G_\fq$
and $I_\fq$ are the decomposition and inertia subgroups of $G_K$ 
at $\fq$.

\section{Level Lowering}

We need a level lowering
result that plays the r\^ole of the Ribet step \cite{RibetLL}
in the proof of Fermat's Last Theorem. 
The following theorem is deduced in \cite{Fermat}
from the work of 
Fujiwara \cite{Fuj},
Jarvis \cite{Jarv} and 
Rajaei \cite{Raj}. 

\begin{thm}\label{thm:levell} 
Let $K$ be a real quadratic field, 
and $E/K$ an elliptic curve of conductor $\cN$.
Let $p$ be a rational prime.
For a prime ideal $\fq$ of $K$ denote by $\Delta_\fq$ the discriminant of a local
minimal model for $E$ at $\fq$.
Let
\begin{equation}\label{eqn:Np}
\cM_p := \prod_{
\substack{\fq \Vert \cN,\\ p \mid \ord_\fq(\Delta_\fq)}
} {\fq}, \qquad\quad \cN_p:=\frac{\cN}{\cM_p} \, .
\end{equation}
Suppose the following:
\begin{enumerate}
\item[(i)] either $p\ge 5$, or $K=\Q(\sqrt{5})$ and $p\ge 7$,
\item[(ii)] $E$ is modular,
\item[(iii)] $\overline{\rho}_{E,p}$ is irreducible,
\item[(iv)] $E$ is semistable at all $\fq \mid p$,
\item[(v)]  $p \mid \ord_\fq(\Delta_\fq)$ for all $\fq \mid p$. 
\end{enumerate}
Then, there is a Hilbert eigenform $\ff$ 
of parallel weight $2$ that is new at level $\cN_p$ and some prime $\varpi$ of $\Q_\ff$
such that $\varpi \mid p$
and $\overline{\rho}_{E,p} \sim \overline{\rho}_{\ff,\varpi}$.
\end{thm}

\section{Scaling and the Odd Part of the Level}\label{sec:odd}
Let $(a,b,c)$ be a non-trivial solution to the
Fermat equation~\eqref{eqn:Fermat}. 
Let $\cG_{a,b,c}$ be as given in \eqref{eqn:cG},
which we think of as the greatest
common divisor of $a$, $b$, $c$. 
An odd prime not dividing $\cG_{a,b,c}$ is a prime
of good or multiplicative reduction for
$E_{a,b,c}$ and does not
appear in the final level $\cN_p$, as we see in due course.
An odd prime dividing
$\cG_{a,b,c}$ exactly once is an additive prime, and 
does appear in $\cN_p$.
To control $\cN_p$, we need to
control $\cG_{a,b,c}$. 

\subsection{Scaling}
We refer to Section~\ref{sec:conv}
for the notation. 
\begin{lem}\label{lem:gcd}
Let $(a,b,c)$ be a non-trivial solution to \eqref{eqn:Fermat}.
There is a non-trivial integral solution 
$(a^\prime,b^\prime,c^\prime)$ to \eqref{eqn:Fermat} 
and some $\fp \in \cH$ such 
that the following hold.
\begin{enumerate}
\item[(i)] For some $\xi \in K^*$, we have
$a^\prime= \xi a$, $b^\prime= \xi b$, $c^\prime=\xi c$.
\item[(ii)] $\cG_{a^\prime,b^\prime,c^\prime}=\fp \cdot \fr^2$ 
where $\fr$ is an odd prime ideal $\ne \fp$.
\item[(iii)] $[a^\prime,b^\prime,c^\prime]=[a,b,c]$.
\end{enumerate}
\end{lem}
\begin{proof}
Recall that $\cH=\{\fp_1,\dotsc,\fp_r\}$
and that $[\fp_1],\dots,[\fp_r]$ represent the cosets
of $\Cl(K)/\Cl(K)^2$. Thus for some $\fp \in \cH$ 
we have $[a,b,c]=[\fp]\cdot [\gb]^2$,
where $\gb$ is a fractional ideal. Now every ideal
class is represented by infinitely many prime ideals.
Thus there is an odd prime ideal $\fr \ne \fp$
such that $[a,b,c]=[\fp]\cdot [\fr]^2$. It follows that
$\fp \cdot \fr^2=(\xi) \cdot \mathcal{G}_{a,b,c}$ for some $\xi \in K^*$. 
Let $a^\prime$, $b^\prime$,
$c^\prime$ be as in (i).  
Note 
\[
(a^\prime)=(\xi) \cdot (a)=\fp \cdot \fr^2 \cdot \cG_{a,b,c}^{-1} \cdot (a)
\]
which is an integral ideal, since $\cG_{a,b,c}$ (by its
definition) divides $a$.
Thus $a^\prime$ is in $\OO_K$ and similarly so are $b^\prime$ and $c^\prime$.
For (ii) and (iii), note that
\[
\cG_{a^\prime,b^\prime,c^\prime}
=a^\prime \OO_K+b^\prime \OO_K+c^\prime \OO_K
=(\xi)\cdot (a \OO_K+b \OO_K+c \OO_K)
=(\xi)\cdot \cG_{a,b,c}=\fp \cdot \fr^2.
\]
\end{proof}

\subsection{Behaviour at odd primes} 
For $u$, $v$, $w \in \OO_K$ such that $uvw\ne 0$ and $u+v+w=0$,
let
\begin{equation}\label{eqn:2tors}
E: y^2=x(x-u)(x+v).
\end{equation}
The invariants $c_4$, $c_6$, $\Delta$, $j$ have their usual meanings 
and are given by:
\begin{equation}\label{eqn:inv}
\begin{gathered}
c_4=16(u^2-vw)=16(v^2-wu)=16(w^2-uv),\\
c_6=-32(u-v)(v-w)(w-u), \qquad 
\Delta=16 u^2 v^2 w^2, \qquad j=c_4^3/\Delta \, .
\end{gathered}
\end{equation}

The following elementary lemma is 
a straightforward consequence of the properties of elliptic
curves over local fields (e.g.\ \cite[Sections VII.1 and VII.5]{SilvermanI}).
\begin{lem}\label{lem:elem}
With the above notation, let $\fq \nmid 2$ be a prime and let
\[
s=\min\{\ord_\fq(u),\ord_\fq(v),\ord_\fq(w) \}.
\]
Write $E_{\mathrm{min}}$ for a local minimal model at $\fq$.
\begin{enumerate}
\item[(i)] $E_{\mathrm{min}}$ has good reduction at $\fq$ if and only if
$s$ is even and 
\begin{equation}\label{eqn:uvweq}
\ord_\fq(u)=\ord_\fq(v)=\ord_\fq(w).
\end{equation}
\item[(ii)] $E_{\mathrm{min}}$ has multiplicative reduction at $\fq$
if and only if $s$ is even and \eqref{eqn:uvweq} fails to hold. In this
case the minimal discriminant 
$\Delta_{\fq}$ at $\fq$ satisfies 
\[
\ord_\fq(\Delta_{\fq})=2 \ord_\fq(u)+2\ord_\fq(v)+2 \ord_\fq(w)-6s.
\]
\item[(iii)] $E_{\mathrm{min}}$ has additive potentially multiplicative reduction at $\fq$ if and only if $s$ is odd and \eqref{eqn:uvweq} fails to hold.
\item[(iv)] $E_{\mathrm{min}}$ has additive potentially good reduction at $\fq$ if and only if $s$ is odd and \eqref{eqn:uvweq} holds.
Moreover, $E$ acquires good reduction over a quadratic extension of $K_\fq$.
\end{enumerate}
\end{lem}

\subsection{The odd part of the level}
Let 
$(a,b,c)$ be a non-trivial solution to the 
Fermat equation \eqref{eqn:Fermat} with odd prime exponent $p$.
Write $E$
for the Frey curve in \eqref{eqn:Frey}. 
Let $\cN$ be the conductor of $E$
and let $\cN_p$
be as defined in \eqref{eqn:Np}. We define the \textbf{even parts} of
$\cN$ and $\cN_p$ by
\[
\cN^\textrm{even}=\prod_{\mP \in S} \mP^{\ord_\mP(\cN)},
\qquad
\cN_p^\textrm{even}=\prod_{\mP \in S} \mP^{\ord_\mP(\cN_p)}.
\]
We define the \textbf{odd parts} of $\cN$ and $\cN_p$ by
\[
\cN^\textrm{odd}=\frac{\cN}{\cN^\textrm{even}},
\qquad
\cN_p^\textrm{odd}=\frac{\cN}{\cN_p^\textrm{even}}.
\]

\begin{lem}\label{lem:cond}
Let  $(a,b,c)$ be a non-trivial solution to the 
Fermat equation \eqref{eqn:Fermat} with odd prime exponent $p$
satisfying $\cG_{a,b,c} = \fp \cdot \fr^2$, where $\fp \in \cH$,
and $\fr$ is an odd prime ideal such that $\fr \ne \fp$. 
Write $E$
for the Frey curve in \eqref{eqn:Frey}. 
Then at all $\fq \notin S \cup \{\fp\}$, 
the local minimal model $E_\fq$ is semistable, and satisfies 
$p \mid \ord_\fq(\Delta_\fq)$.
Moreover,
\begin{equation}\label{eqn:cnp}
\cN^{\textrm{odd}}=
\fp^{2}\cdot
\fr^{\text{$0$ or $1$}} 
\cdot 
\prod_{\substack{\fq \mid abc \\ 
\fq \notin S \cup \{\fp, \fr\} 
}} \fq, 
\qquad 
\qquad
\cN_p^{\textrm{odd}}= 
\fp^{2} .
\end{equation}
\end{lem}
\begin{proof}
Clearly, if $\fq \nmid 2abc$ then $E$ has good reduction at $\fq$, hence $\fq \nmid \cN$, $\cN_p$.  Note also that
\[
\min\{ \ord_\fr(a^p), \ord_\fr(b^p), \ord_\fr(c^p)\}=2p.
\]
By Lemma~\ref{lem:elem}, $E$ has good or multiplicative
reduction at $\fr$, and in either case $p \mid \ord_\fr(\Delta_\fr)$,
proving also the correctness of the exponents of $\fr$
in $\cN$ and $\cN_p$.

Recall 
that $\fp \in \cH$ satisfy $\fp=1\cdot \OO_K$, or  $\fp$
is an odd prime ideal. In the former case there is nothing
to prove, so suppose that $\fp$ is an odd prime ideal.
As $E$ has full $2$-torsion over $K$, the wild part of the conductor
of $E/K$ at $\fp$ vanishes (see \cite[page 380]{SilvermanII}).
Moreover,
\[
\min\{ \ord_\fp(a^p), \ord_\fp(b^p), \ord_\fp(c^p)\}=p.
\]
By Lemma~\ref{lem:elem}, $E/K$ has additive reduction
at $\fp$. Thus the exponent of $\fp$ in both 
$\cN$ and $\cN_p$ is $2$.

Suppose that $\fq \mid abc$ and $\fq \notin S \cup \{\fp,\fr\}$.
Since $\cG_{a,b,c}= \fp \cdot \fr^2$, 
the prime $\fq$ divides precisely one of $a$, $b$, $c$.
From \eqref{eqn:inv}, $\fq \nmid c_4$ so the model \eqref{eqn:Frey}
is minimal and
has multiplicative reduction at $\fq$, and $p \mid \ord_\fq(\Delta)$.
By \eqref{eqn:Np}, we see that $\fq \nmid \cN_p$.
\end{proof}


\begin{cor}\label{cor:condodd}
Let $2 \le d \le 23$ be squarefree and let $K=\Q(\sqrt{d})$.
Let 
$(a,b,c)$ be a non-trivial solution to the 
Fermat equation \eqref{eqn:Fermat}.
We may scale $(a,b,c)$ so that it remains integral, and 
\[
\cG_{a,b,c} = \fp \cdot \fr^2, \qquad \cN_p^\textrm{odd}=\fp^2
\] 
where 
\begin{enumerate}
\item[(a)] if $d \ne 10$, $15$ then $\fp=1 \cdot \OO_K$;
\item[(b)] if $d =10$ then $\fp=1 \cdot \OO_K$ or
$\fp=(3,1+\sqrt{10})$;
\item[(c)] if $d=15$ then $\fp=1 \cdot \OO_K$ or
$\fp=(3,\sqrt{15})$;
\end{enumerate}
and $\fr$ is an odd prime ideal such that $\fr \ne \fp$. 
\end{cor}
\begin{proof}
The corollary follows immediately from Lemmas~\ref{lem:gcd},~\ref{lem:cond}
and a computation of a suitable set $\cH$
for each of the quadratic fields.
\end{proof}

\section{Scaling by units and the even part of the level}\label{sec:scaling}
In the previous section we scaled $(a,b,c)$
so that $\cG_{a,b,c}=\fp \cdot \fr^2$, where $\fp \in \cH$,
and we computed the odd parts of the conductor $\cN$
and level $\cN_p$. 
Let $\OO_K^*$ be the unit group of $K$.
In this section we study the effect 
on $\cN$ and $\cN_p$
of
scaling $(a,b,c)$ by units.
Note that scaling $(a,b,c)$ 
by units does not affect $\cG_{a,b,c}$; it is plain from
the proofs in the previous section that this
leaves the odd parts of $\cN$, $\cN_p$ unchanged.
Applying an even permutation to $(a,b,c)$ results in
an isomorphic Frey curve, whereas applying an odd permutation
replaces the Frey curve by its twist by $-1$, and so has
the same effect as scaling $(a,b,c)$ by $-1$.

\begin{lem}\label{lem:inert}
Suppose $K$ is a quadratic field and $2$ is inert in $\OO_K$.
Let $\mP=2 \OO_K$, 
and suppose $\mP \nmid abc$.
Then after suitably permuting $(a,b,c)$ we have
$\ord_{\mP}(\cN)=4$. Moreover, $E$ has potentially good reduction at $\mP$.
\end{lem}
\begin{proof}
We can write $K=\Q(\sqrt{d})$ where $d \equiv -3 \pmod{8}$.
Thus $\OO_\mP=\Z_2[\omega]$ where $\omega^2+\omega+1=0$.
The residue field of $\mP$ is $\F_2[\tilde{\omega}] \cong \F_4$.
Write $A=a^p$, $B=b^p$, $C=c^p$. 
Now $\mP \nmid ABC$ and $A+B+C=0$. 
Then $A$, $B$, $C$ are congruent modulo $\mP$
to $1$, $\omega$, $\omega^2$ in some order. By rearranging,
we may suppose
that $C \equiv \omega^2 \pmod{\mP}$, and we will decide later
on which of $A$, $B$ are congruent to $1$ and $\omega$ modulo $\mP$.
Let $E$ denote the Frey curve in \eqref{eqn:Frey}. 
It follows from \eqref{eqn:inv} that $\ord_{\mP}(\Delta)=4$
and $\ord_\mP(c_4) =4$. In particular, $\ord_\mP(j)=8$,
and so $E$ has potentially good reduction at $\mP$.
Furthermore, the Frey curve is minimal at $\mP$ and has additive reduction.
We will follow
the steps of Tate's algorithm as in \cite[page 366]{SilvermanII}.
Let $\tilde{E}$ denote the reduction of $E$ modulo $\mP$.
It is easy to check that the point $(\tilde{C},\tilde{1})$
is singular on $\tilde{E}$. Now we shift the model $E$,
replacing $X$ by $X+C$ and $Y$ by $Y+1$, which has the effect
of sending sending the point $(\tilde{C},\tilde{1})$ on
the special fibre to $(\tilde{0},\tilde{0})$. Write
$a_1,\dotsc,a_6$ for the $a$-invariants of the resulting model.
Then
\[
a_6=C^3+(B-A) C^2-ABC-1.
\]
By Step 3 of Tate's Algorithm we know that if $\mP^2 \nmid a_6$
then $\ord_\mP(\cN)=\ord_\mP(\Delta)=4$. Suppose $\mP^2 \mid a_6$.
Now swapping $A$, $B$, replaces $a_6$ by
\[
a_6^\prime=C^3+(A-B) C^2-ABC-1.
\]
Observe that $\ord_\mP(a_6^\prime-a_6)=\ord_\mP(2(A-B)C^2)=1$.
Hence $\mP^2 \nmid a_6^\prime$. Thus we may always permute $A$, $B$, $C$
so that $\ord_\mP(\cN)=4$. 
\end{proof}
\noindent \textbf{Remark.} Under the hypotheses of Lemma~\ref{lem:inert},
it follows from Ogg's formula \cite[Section IV.11]{SilvermanII} that
the possible exponents of $\mP$ in the conductor are $2$, $3$, $4$.
Lemma~\ref{lem:inert} shows that we can always permute the solution so
that the exponent of $\mP$ in the conductor is $4$, avoiding having to compute
newforms at the smaller levels. Unfortunately, we have found that it is not
possible by permuting the solution to force the exponent to be smaller in general.

\begin{lem}\label{lem:mult2}
Suppose $K$ is a quadratic field and let $\mP \in S$. Suppose $(a,b,c)$ is a
non-trivial solution to the Fermat equation, with $\mP \nmid \cG_{a,b,c}$.  
The Frey curve $E$ has potentially multiplicative reduction
at $\mP$ if and only if
\begin{enumerate}
\item[(a)] either $f(\mP/2)=1$ (i.e. $2$ splits or ramifies in $K$),
\item[(b)] or $f(\mP/2)=2$ (i.e. $2$ is inert in $K$)
and $\mP \mid abc$.
\end{enumerate}
Moreover, if the reduction at $\mP$ is multiplicative then
$p \nmid \ord_\mP(\Delta_\mP)$.
\end{lem}
\begin{proof}
Suppose (a) or (b) holds. We claim that $\mP \mid abc$.
If (b) holds this is true by hypothesis. If (a), then
the residue
field at $\mP$ is $\F_2$. It follows from
$a^p+b^p+c^p=0$ that $\mP$ divides at least one of $a$, $b$, $c$, establishing our claim.
Moreover, as 
$\mP \nmid \cG_{a,b,c}$, we see that $\mP$ divides precisely one of $a$, $b$, $c$.
Let $t=\ord_\mP(abc) \ge 1$. By \eqref{eqn:inv},
\[
\ord_\mP(c_4)=4 \ord_\mP(2),
\qquad
\ord_\mP(\Delta)=4 \ord_\mP(2)+2p t.
\]
Thus
\begin{equation}\label{eqn:jval}
\ord_\mP(j)=8 \ord_\mP(2)-6 p t <0
\end{equation}
as $p \ge 3$. Thus we have potentially multiplicative reduction at $\mP$.
The converse follows from Lemma~\ref{lem:inert}.

To complete the proof suppose that the reduction is multiplicative, and let
$c_4^\prime$, $c_6^\prime$ and $\Delta^\prime=\Delta_\mP$ be
the corresponding invariants of a local minimal model. Now
$\mP \nmid c_4^\prime$, but $j^\prime={c_4^\prime}^3/\Delta^\prime$.
From \eqref{eqn:jval}, $p \nmid \ord_\mP(j)$, and hence $p \nmid \ord_\mP(\Delta_\mP)$.
\end{proof}

\begin{lem}\label{lem:integration}
Let $K_\mP$ be a local field, and $E$ and elliptic curve 
over $K_\mP$ with potentially multiplicative reduction.
Let
$c_4$, $c_6$ be the usual $c$-invariants of $E$.
Let $L=K_\mP(\sqrt{-c_6/c_4})$ and let $\Delta(L/K_\mP)$
be the discriminant of this local extension.
Then the conductor of $E/K_\mP$ is
\[
f(E/K_\mP)=\begin{cases}
1 & \text{if $\ord_\mP(\Delta(L/K_\mP))=0$} \\
2 \ord_\mP(\Delta(L/K_\mP)) & \text{otherwise}.
\end{cases}
\]
\end{lem}
\begin{proof}
Let $E^\prime$ be the quadratic twist of $E$
by $-c_6/c_4$. Then $E^\prime$ is a Tate curve (see
for example \cite[Section V.5]{SilvermanII}).
The Lemma now follows from \cite[Section 18]{Rohrlich2}.
\end{proof}

\begin{lem}
Let $(a,b,c)$ be a non-trivial solution to the Fermat equation,
such that $\cG_{a,b,c}$ is odd. 
Suppose $2$ is either split or ramified in $K$, or that
$2$ is inert and $2 \mid abc$. Let
\[
\gb=\prod_{\mP \in S} \mP^{2\ord_\mP(2)+1},
\]
and write
\[
\Phi : \OO_K^* \rightarrow (\OO_K/\gb)^*/((\OO_K/\gb)^*)^2
\]
for the natural map.
Choose a set  
$\lambda_1, \dotsc, \lambda_k \in \OO_K\backslash \gb$ 
that represent the elements of the cokernel of 
$\Phi$. For $1 \le i \le k$, and for $\mP \in S$, let
$\Delta_\mP^{(i)}$ be the discriminant of the local
extension $K_\mP(\sqrt{\lambda_i})/K_\mP$,
and let 
\[
\epsilon_\mP^{(i)}
=\begin{cases}
1 & \text{if $\ord_\mP(\Delta_\mP^{(i)})=0$}\\
2 \ord_\mP(\Delta_\mP^{(i)}) & \text{otherwise}.
\end{cases}
\]
Then we may scale $(a,b,c)$ by an element of $\OO_K^*$
so that for some $i$, and for every $\mP \in S$,
we have $\ord_\mP(\cN)=\epsilon_\mP^{(i)}$.
\end{lem}
\begin{proof}
Write $\OO = \OO_K$. By Lemma~\ref{lem:mult2} 
we have potentially multiplicative reduction all the primes $\mP \in S$.
Write $c_4$, $c_6$ for the usual $c$-invariants of the Frey curve $E$.
Since $\cG_{a,b,c}$ is odd, but $\mP \mid abc$ for all $\mP \in S$,
we have from \eqref{eqn:inv} and the relation $a^p+b^p+c^p=0$
that $\ord_\mP(c_4)=4 \ord_\mP(2)$ and $\ord_\mP(c_6)=6 \ord_\mP(2)$.
Write $\gamma=-c_6/4 c_4$. Then $\gamma \in \OO_\mP^*$ for all $\mP \in S$,
and $K_\mP(\sqrt{\gamma})=K_\mP(\sqrt{-c_6/c_4})$.
Now 
the exponent of $\mP$ in the conductor $\cN$ of the Frey curve
can be expressed by Lemma~\ref{lem:integration} in terms of the
discriminant of the extension $K_\mP(\sqrt{\gamma})/K_\mP$.

We shall make use of the isomorphism
\[
(\OO/\gb)^*/((\OO/\gb)^*)^2
\cong
\prod_{\mP \in S} \OO_\mP^*/(\OO_\mP^*)^2
\]
which follows from the Chinese Remainder Theorem, and Hensel's Lemma.
Observe that
scaling $(a,b,c)$ by a unit $\eta \in \OO_K^*$ scales $\gamma$
by $\eta^p$. Now, as $p$ is odd, it follows from
the definition of $\Phi$ and the above isomorphism that we can scale
$(a,b,c)$ by some $\eta \in \OO_K^*$
so that there is some
$1 \le i \le k$ with $\gamma/\lambda_\mP^{i}$ a square
in $\OO_\mP$ for each $\mP \in S$. Therefore, 
\[
 K_\mP(\sqrt{\gamma}) = K_\mP\left(\sqrt{\lambda_\mP^{i}}\right)
\]
and the lemma follows from  Lemma~\ref{lem:integration}. 
\end{proof}

\noindent \textbf{Remark}. 
Let $u$ be the fundamental unit of the real quadratic field $K$.
Observe that if $\lambda \in \OO_K\backslash \gb$
represents an element of the cokernel of $\Phi$,
then for every integer $k$, 
the same element of the cokernel is
also represented by $\lambda^\prime=\pm u^k \lambda$.
The local extension $K_\mP(\sqrt{\lambda^\prime})/K_\mP$
depends only on the choice of sign $\pm$ and the parity 
of $k$. To keep the even part of the level small, we
replace each representative $\lambda$ by whichever one  of $\lambda$, $-\lambda$, $u \lambda$, $-u \lambda$ that minimizes the norm
of the even part of the level $\cN_p$.

\section{Possibilities for $\cN_p$}

\begin{table}
\begin{tabular}{||c|c|c|c||}
\hline\hline
\multirow{2}{*}{$d$}
&
\multirow{2}{*}{$S$} 
&
\multirow{2}{*}{$\lambda$s}
&
$\cN^{\textrm{even}}=$
\\
&
&
&
$\cN_p^{\textrm{even}}$
\\
\hline\hline
$2$
&
$\mP=(\sqrt{2})$
&
$1, \quad -1-2\sqrt{2}$
& 
$\mP$
\\
\hline\hline
\multirow{2}{*}{$3$}
&
\multirow{2}{*}{$\mP=(1+\sqrt{3})$}
&
$1$ & $\mP$ 
\\ \cline{3-4}
& & $-1+2\sqrt{3}$ 
& $\mP^4$
\\
\hline\hline
\multirow{2}{*}{$5$}
&
\multirow{2}{*}{$\mP=(2)$}
&
$1, \quad -5+2\sqrt{5}$
&
$\mP$
\\ \cline{3-4}
&
&
$\mP \nmid abc$
&
$\mP^4$
\\
\hline\hline
\multirow{2}{*}{$6$}
&
\multirow{2}{*}{$\mP=(-2+\sqrt{6})$}
&
$1$
& 
$\mP$ \\ \cline{3-4}
& &
$1+\sqrt{6}$
&
$\mP^8$
\\
\hline\hline
\multirow{2}{*}{$7$}
&
\multirow{2}{*}{$\mP=(3+\sqrt{7})$}
&
$1, \quad 21-8\sqrt{7}$ 
&
$\mP$
\\ \cline{3-4}
&
&
$-1+2\sqrt{7}, \quad -5+2\sqrt{7}$
&
$\mP^4$
\\
\hline
\hline
$10$
&
$\mP=(2,\sqrt{10})$
&
$1, \qquad 7-2\sqrt{10}$
&
$\mP$
\\
\hline
\hline
\multirow{2}{*}{$11$}
&
\multirow{2}{*}{$\mP=(3+\sqrt{11})$}
&
$1$
&
$\mP$\\ \cline{3-4}
& &
$-1+2\sqrt{11}$
&
$\mP^4$
\\
\hline\hline
\multirow{2}{*}{$13$}
&
\multirow{2}{*}{$\mP=(2)$}
&
$1,\quad -5+2\sqrt{13}$
&
$\mP$
\\ \cline{3-4}
&
&
$\mP \nmid abc$
&
$\mP^4$
\\
\hline\hline
\multirow{2}{*}{$14$}
&
\multirow{2}{*}{$\mP=(4+\sqrt{14})$}
&
$1,\qquad -3$
& 
$\mP$
\\ \cline{3-4}
&
&
$1+\sqrt{14},\qquad -3+\sqrt{14}$
&
$\mP^8$
\\
\hline\hline
\multirow{2}{*}{$15$}
&
\multirow{2}{*}{$\mP=(2,1+\sqrt{15})$}
&
$1, \quad -15+4 \sqrt{15}$
&
$\mP$
\\ \cline{3-4}
&
&
$-1+2\sqrt{15},\quad 7-2\sqrt{15}$
& 
$\mP^4$
\\
\hline\hline
$17$
&
$\mP_1=\left(\frac{3+\sqrt{17}}{2}\right),
\;
\mP_2=\left(\frac{3-\sqrt{17}}{2}\right)$
&
 $1, \; 17-4\sqrt{17}, \; -9+2\sqrt{17}, \; -5+2\sqrt{17}$
&
$\mP_1 \cdot \mP_2$
\\
\hline \hline
\multirow{2}{*}{$19$}
&
\multirow{2}{*}{$\mP=(13+3\sqrt{19})$}
&
$1$ 
&
$\mP$
\\ \cline{3-4}
&
&
$-1+2\sqrt{19}$
&
$\mP^4$
\\
\hline\hline
\multirow{3}{*}{$21$}
&
\multirow{3}{*}{$\mP=(2)$}
&
$1, \quad -5+2\sqrt{21}$
&
$\mP$ 
\\ \cline{3-4}
&
&
$(7-\sqrt{21})/2,\quad (3+3\sqrt{21})/2$
&
$\mP^4$
\\ \cline{3-4}
&
&
$\mP \nmid abc$
&
$\mP^4$
\\
\hline\hline
\multirow{2}{*}{$22$}
&
\multirow{2}{*}{$\mP=(14+3\sqrt{22})$}
&
$1$
& 
$\mP$ 
\\ \cline{3-4}
&
&
$1+\sqrt{22}$
&
$\mP^8$
\\
\hline\hline
\multirow{2}{*}{$23$}
&
\multirow{2}{*}{$\mP=(5+\sqrt{23})$}
&
$1,\quad
   115+ 24 \sqrt{23}$
&
$\mP$
\\ \cline{3-4}
&
&
$-1+2\sqrt{23}, \quad -163-34 \sqrt{23}$
&
$\mP^4$
\\
\hline\hline
\end{tabular}
\caption{Quantities required for Corollary~\ref{cor:cond} and its proof.}
\label{table:Neven}
\end{table}

\begin{cor}\label{cor:cond}
Let $2 \le d \le 23$ be squarefree and let $K=\Q(\sqrt{d})$.
Let 
$(a,b,c)$ be a non-trivial solution to the 
Fermat equation \eqref{eqn:Fermat} with odd prime exponent $p$.
We may scale $(a,b,c)$ so that it remains integral,  
$\cG_{a,b,c}$ and $\cN_p^{\textrm{odd}}$
are as in Corollary~\ref{cor:condodd} and 
$\cN_p^{\textrm{even}}=\cN^{\textrm{even}}$
is as given in Table~\ref{table:Neven}.
\end{cor}

\section{Irreducibility}\label{sec:irred}

\begin{lem}\label{lem:irred}
Let $p\ge 17$, and let $2 \le d \le 23$ be squarefree. Let $K=\Q(\sqrt{d})$.
Let $(a,b,c)$ be a non-trivial solution to the Fermat equation~\eqref{eqn:Fermat} scaled as in Corollary~\ref{cor:cond}.
Then $\overline{\rho}_{E,p}$ is irreducible.
\end{lem}
\begin{proof}
Suppose that $\overline{\rho}_{E,p}$ is reducible. 
As $E$ has non-trivial $2$-torsion,
it gives rise to a $K$-point on $X_0(2p)$.
The quadratic points on $X_0(34)$ have been determined by Ozman \cite{Ekin}.
These are all defined over $\Q(i)$, $\Q(\sqrt{-2})$ and $\Q(\sqrt{-15})$.
Thus we suppose $p \ge 19$.

We can write
\[
\overline{\rho}_{E,p} \sim 
\begin{pmatrix}
\theta & * \\
0 & \theta^\prime
\end{pmatrix},
\]
where $\theta$ and $\theta^\prime$ are characters $G_K \rightarrow \F_p^*$.
By replacing $E$ by a $p$-isogenous elliptic curve we may swap $\theta$
and $\theta^\prime$ as we please. The characters $\theta$ and $\theta^\prime$ are unramified away from
$p$ and the additive primes for $E$ (see the proof of \cite[Lemma 1]{KrausQuad}).
We shall write $\cN_\theta$ and $\cN_{\theta^\prime}$ for the conductors
of $\theta$ and $\theta^\prime$.
Since $\theta^\prime=\chi/\theta$, where $\chi$ is the cyclotomic
character,
we see that $\theta$ and $\theta^\prime$ have the same conductor
away from $p$. It is easy to show that if $\fq \nmid p$ is a prime of
additive reduction, then $\ord_\fq (\cN)$ is even, and 
\[
\ord_\fq(\cN_\theta)=\ord_\fq(\cN_{\theta^\prime})=\frac{1}{2} \ord_\fq(\cN).
\]
We now divide into two cases according to whether one of $\theta$, $\theta^ \prime$ is unramified at $p$.

\smallskip

\noindent (i) Suppose first that $p$ is coprime to $\cN_\theta$ or $\cN_{\theta^\prime}$. After swapping $\theta$ and $\theta^\prime$
we can assume that $p$ is coprime to $\cN_\theta$. Thus $\cN_\theta$
is the squareroot of the additive part of the conductor $\cN$. From
Lemma~\ref{lem:cond}, we know that the odd additive part is $\fp^2$
where the possibilities for $\fp$ are as in Corollary~\ref{cor:condodd}.
The even additive part of $\cN$ can be deduced from Table~\ref{table:Neven}.
Thus for each $d$ we have a small list of possibilities for
$\cN_\theta$. Let $\infty_1$ and $\infty_2$ be the two real
places from $K$. 
It follows that $\theta$ is a character of the
ray class group
for the modulus $\cN_\theta \infty_1 \infty_2$. Using
{\tt Magma} we computed this ray class group in all cases,
and found it to be one of the following groups
\[
0, \quad \Z/2\Z, \quad \Z/4\Z, \quad 
\Z/2\Z \times \Z/2\Z, \quad
\Z/2\Z \times \Z/4\Z.
\]
The order of $\theta$ divides the exponent of the group, and so
it is $1$, $2$ or $4$. If $\theta$ has order $1$, then $E$
has a point of order $p$ over $K$. The possibilities for $p$-torsion over
quadratic fields have been determined by Kamienny, Kenku and Momose 
(c.f. \cite[Theorem 3.1]{Kamienny}) and their results imply that $p \le 13$ giving a contradiction.
If $\theta$ has order $2$, then $E$ has point of order $p$
over a quadratic extension $L/K$. The field $L$ has degree $4$ over $\Q$.
The possibilities for $p$-torsion over quartic fields have
been determined by Derickx, 
Kamienny, Stein and Stoll \cite{torsion} and their results
imply $p \le 17$, again giving a contradiction.
Suppose $\theta$ has order $4$. Let $L$ be the unique quadratic
extension of $K$ cut out by $\theta^2$. Now $\phi=\theta \vert_{G_L}$ 
is a quadratic character. Twisting $E/L$ by $\phi$ gives an elliptic
curve defined over $L$ with a point of order $p$. As before, $p \le 17$.

\medskip

\noindent (ii) Suppose that neither $\cN_\theta$, $\cN_{\theta^\prime}$
is coprime to $p$. Let $\gp \mid p$ be a prime of $K$. Then $E$
is semistable at $\gp$. If $\gp/p$ is inert, it follows from \cite[Lemma
1]{KrausQuad} that $\gp$ divides at most one of $\cN_\theta$,
$\cN_{\theta^\prime}$ and so we are in case (i). Thus we may assume that $p$
splits or ramifies in $K$. 

(iia) Suppose now that $p$ ramifies in $K$. 
This means that $d=p$ and $d$ is either $19$ or $23$. 
Let $\gp$ be the unique prime above $p$ in $K$. 
If $E$ has good ordinary or multiplicative reduction at $\gp$, 
then $\gp$ divides at most one of $\cN_\theta$, $\cN_{\theta^\prime}$ 
and we are finished by case (i). Thus suppose $E$ has good 
supersingular reduction at $\gp$.
We will now apply Proposition~\ref{pp:supersing} below to show 
this cannot happen. The field $K=\Q(\sqrt{d})$ has a prime 
$\mP \mid 2$ with residue field $\F_2$, 
and so by Lemma~\ref{lem:mult2}
 this is a prime of potentially multiplicative reduction for
$E$. In the notation of 
Proposition~\ref{pp:supersing},
$\sqrt{\mathcal{A}}=(1)$ or $\mP^2$. In all cases,
the ray class group of modulus
$\sqrt{\mathcal{A}} \infty_1 \infty_2$ is $\Z/2\Z$,
and the proposition implies that
$4=\norm(\mP)^2 \equiv 1 \pmod{p}$.  As 
$p=19$ or $23$, we have a contradiction and so $E$ cannot be supersingular.

(iib) Suppose $p$ splits as $\gp \gp^\prime$. Moreover, we can also assume $\gp \mid \cN_\theta$, $\gp \nmid \cN_{\theta^\prime}$ and $\gp^\prime \nmid \cN_{\theta}$, $\gp^\prime \mid \cN_{\theta^\prime}$.
The primes $\gp$, $\gp^\prime$ are unramified. It follows \cite[Proposition 12]{Serre72} that $E$
has either good ordinary reduction or multiplicative reduction at these.
Thus
$\theta \vert_{I_\gp}=\chi \vert_{I_\gp}$ and  
$\theta^\prime \vert_{I_{\gp^\prime}}=\chi \vert_{I_{\gp^\prime}}$. 
We shall write down a small integer $n>0$ such that
$\theta^n$ is unramified away from $\gp$. If $\fq \nmid p$ is
a prime of potentially multiplicative reduction then
$\theta^2$ is unramified at $\fq$. Furthermore, for our Frey curve $E$,
the only odd additive prime is $\fq=\fp$, and
part (iv) of Lemma~\ref{lem:elem}
implies that 
$\#\overline{\rho}_{E,p}(I_\fq)=2$, and so
$\theta^2$ is unramified at $\fq$. 
We are left with primes $\fq \mid 2$
of potentially good reduction. These only arise for $d=5$, $13$, $21$,
and in these cases $\ord_\fq(\Delta_\fq)=4$. It follows from \cite[Theorem 3]{kraus2}
that $\overline{\rho}_{E,p}(I_\fq)$ is either cyclic of order $3$, $6$
or isomorphic to $\SL_2(\F_3)$ and so has order $24$.
The last case cannot occur, as $\SL_2(\F_3)$ is non-abelian,
and any non-abelian reducible subgroup of $\GL_2(\F_p)$
has an element of order $p$. It follows that $\theta^6$
is unramified at $\fq$. Letting $n=6$ for $d=5$, $13$, $21$,
and $n=2$ for other values of $d$, we conclude that the
character $\theta^n$ is unramified away from $\gp$,
 and that $\theta^n \vert_{I_\gp}
=\chi^n \vert_{I_\gp}$. Let $u$ be a generator of the subgroup
of totally positive units in $\OO_K^*$. It follows (c.f.\ \cite[page 249]{KrausQuad})
that $p \mid \norm(u^n-1)$. We computed the factorization of $\norm(u^n-1)$
for our values of $d$, and found that none are divisible by primes $p\ge 19$,
except when  $d=p=19$ or $d=p=23$. However, in 
these cases $p$ ramifies in the field $K=\Q(\sqrt{d})$, and so are covered by case (iia).
\end{proof}

\begin{prop} \label{pp:supersing} 
Let $d=p \geq 5$ be a prime and $\gp$ be the unique prime in $K(\sqrt{d})$
above $p$. Let $E/K$ be an elliptic curve and denote by $\mathcal{A}$ the
additive part of its conductor. Suppose that $E$ has good supersingular
reduction at $\gp$ and potentially multiplicative reduction at some prime $\fq_0
\neq \gp$. Suppose further that $\bar{\rho}_{E,p}$ is reducible. Therefore,
$\mathcal{A}$ is a square and we let $n$ be the exponent of the ray class group
for the modulus $\sqrt{\mathcal{A}} \infty_1 \infty_2$. Then, $\norm(\fq_0)^{n}
\equiv 1 \pmod{p}$.  
\end{prop}
\begin{proof} 
Suppose $\bar{\rho}_{E,p}$ is reducible and let $\theta, \theta^\prime$ be as
in the proof of Lemma~\ref{lem:irred}.  Since $E$ is supersingular at $\gp$, it
is well known \cite[Proposition 10]{Serre72}  that 
\[
\overline{\rho}_{E,p} \vert_{I_\gp}
\sim 
\begin{pmatrix}
\psi_2^2 & 0\\
0 & \psi_2^{2p}
\end{pmatrix},
\qquad
\text{or}
\qquad
\overline{\rho}_{E,p} \vert_{I_\gp}
\sim 
\begin{pmatrix}
\psi_1 & 0\\
0 & \psi_1
\end{pmatrix},
\]
where $\psi_m$ denotes a degree $m$ fundamental character.
The character $\psi_2^2$ is not $\F_p$-valued, which gives
a contradiction with reducibility. We thus suppose that inertia
acts via the level 1 fundamental characters.

Clearly $\epsilon=\theta/\theta^\prime$ is unramified
at $\gp$. Moreover, as $\epsilon=\theta^2/\chi$ where $\chi$ is the cyclotomic
character, it follows that for $\fq \nmid p$ an additive prime,
\[
\ord_\fq(\cN_\epsilon) \le \ord_\fq(\cN_{\theta}) =  \frac{1}{2}
\ord_\fq(\mathcal{A}).  \]
Therefore, the exponent of the ray class group of modulus $\cN_\epsilon \infty_1
\infty_2$ is a divisor of $n$. Thus $\epsilon^n=1$. Let $\sigma_{\fq_0}$ be the
Frobenius element of $G_K$ at $\fq_0$. Since $\fq_0$ is of potentially
multiplicative reduction the possible pairs of eigenvalues of
$\overline{\rho}_{E,p}(\sigma_{\fq_0})$ are $(1, \norm(\fq_0))$ or $(-1,
-\norm(\fq_0))$ and they correspond to the values
of $\theta(\sigma_{\fq_0})$ and $\theta^\prime(\sigma_{\fq_0})$ up to reordering.
Thus, 
\[
1=\epsilon^n(\sigma_{\fq_0})=\theta(\sigma_{\fq_0})^n/\theta^\prime(\sigma_{\fq_0})^n \equiv
 \norm(\fq_0)^{\pm n}
\pmod{p} .
\]
\end{proof}

\section{Proof of Theorem~\ref{thm:main}}\label{sec:proof1}
For now let $2 \le d \le 23$ be squarefree, and let $K=\Q(\sqrt{d})$.
We would like to show that the equation $x^n+y^n=z^n$ has
only trivial solutions in $K$ for $n \ge 4$,
although as we will see in due course, our proof strategy 
fails for $d=5$ and $d=17$. As in the introduction, we reduce
to showing that the Fermat equation~\eqref{eqn:Fermat}
has no non-trivial solutions $(a,b,c)$ in $\OO_K$
with prime exponent $p \ge 17$. Now suppose $(a,b,c)$ 
is a non-trivial solution with $p \ge 17$, and
scale this as in
Corollary~\ref{cor:cond}. Let $E=E_{a,b,c}$ be
the Frey curve given by~\eqref{eqn:Frey}, and let
$\overline{\rho}_{E,p}$ be its mod $p$ representation. We know from
Lemma~\ref{lem:irred} that $\overline{\rho}_{E,p}$ is irreducible.
We now apply Theorem~\ref{thm:levell} to deduce that there
is a cuspidal Hilbert newform $\ff$ over $K$ of 
weight $(2,2)$ and level $\cN_p$ (one of the levels predicted
by Corollary~\ref{cor:cond}) such that 
$\overline{\rho}_{E,p} \sim \overline{\rho}_{\ff,\varpi}$
for some prime $\varpi \mid p$ of $\Q_\ff$.

\begin{lem}
Let $\fq \nmid \cN_p$ be a prime of $K$, and let 
\[
\cA=\{ a \in \Z 
\quad
: 
\quad \lvert a \rvert \le 2 \sqrt{\norm(\fq)}, \qquad 
\norm(\fq)+1-a \equiv 0 \pmod{4}\}.
\]
If 
$\overline{\rho}_{E,p} \sim \overline{\rho}_{\ff,\varpi}$
then $\varpi$ divides the principal ideal
\[
B_{\ff,\fq}=\norm(\fq) ((\norm(\fq)+1)^2-a_\fq(\ff)^2) 
\prod_{a \in \cA} (a- a_\fq(\ff)) \cdot \OO_{\Q_\ff} \, .
\]
\end{lem}
\begin{proof}
If $\fq \mid p$, then $\norm(\fq)$ is a power of $p$,
and so $\varpi \mid p$ divides $B_{\ff,\fq}$. Thus we may
suppose $\fq \nmid p$. By assumption $\fq \nmid \cN_p$. From the
definition of $\cN_p$ in \eqref{eqn:Np}, the prime
$\fq$
is of good or multiplicative reduction for $E$. 
If $\fq$ is a prime of good reduction for $E$, then
$a_\fq(E) \equiv a_\fq(\ff) \pmod{\varpi}$. By the Hasse--Weil
bounds we know that $\lvert a_\fq(E) \rvert \le 2 \sqrt{\norm(\fq)}$.
Moreover, as $E$ has full $2$-torsion (and $\fq \nmid 2$
as $\fq \nmid \cN_p$) we have $4 \mid \# E(\F_\fq)$.
Thus $a_\fq(E) \in \cA$ and so $\varpi \mid B_{\ff,\fq}$.
Finally, suppose $\fq$ is a prime of multiplicative
reduction for $\cN_p$. Then, comparing the traces of
the images of Frobenius at $\fq$ under $\overline{\rho}_{E,p}$
and $\overline{\rho}_{\ff,\varpi}$ we have
\[
\pm (\norm(\fq)+1) \equiv a_\fq(\ff) \pmod{\varpi}.
\] 
It follows that $\varpi$ divides $B_{\ff,\fq}$ in this case too.
\end{proof}

Using {\tt Magma} we computed the newforms $\ff$
at the predicted levels, the fields $\Q_\ff$,
and eigenvalues $a_\fq(\ff)$
at primes $\fq$ of $K$ small norm. 
We computed for each $\ff$
at level $\cN_p$ the ideal
\[
B_\ff:= \sum_{\fq \in T} B_{\ff,\fq}
\]
where $T$ is the set of prime ideals $\fq \nmid \cN_p$ of $K$ 
with norm $<60$. Let $C_\ff:=\norm_{\Q_\ff/\Q}(B_\ff)$.
If $\overline{\rho}_{E,p} \sim \overline{\rho}_{\ff,\varpi}$
then by the above lemma, $\varpi \mid B_\ff$ and so 
$p \mid C_\ff$. Hence, the isomorphism
$\overline{\rho}_{E,p} \sim \overline{\rho}_{\ff,\varpi}$
is impossible if $p \nmid C_\ff$. Thus, the newforms 
satisfying $C_\ff = 0$ are the problematic ones.
We computed $C_\ff$ for all newforms $\ff$
at the predicted levels, and found only three fields where
$C_\ff = 0$ for some $\ff$. All the others produced values of $C_\ff$ that
are not divisible by primes $p \ge 17$. Thus to complete the proof
we have to deal with the cases where $C_\ff=0$; these
are as follows:
\begin{enumerate}
\item[(i)] $K=\Q(\sqrt{3})$, $\cN_p=(1+\sqrt{3})^4$.
Here $\ff$ is the unique newform at level $\cN_p$.
It satisfies $\Q_\ff=\Q$ and corresponds to
the elliptic curve
\[
E^\prime \; : \; y^2=x(x+1)(x+8+4\sqrt{3})
\]
of conductor $(1+\sqrt{3})^4$.
\item[(ii)] $K=\Q(\sqrt{5})$, $\cN_p=(2)^4$.
There are three newforms at level $\cN_p$,
and all there satisfy $\Q_\ff=\Q$. For 
all three newforms $C_\ff=0$.
\item[(iii)] $K=\Q(\sqrt{17})$, $\cN_p=(2)$.
Here $\ff$ is the unique newform at level $\cN_p$.
It satisfies $\Q_\ff=\Q$ and corresponds
to the elliptic curve
\[
W \; : \;
y^2= x(x-4+\sqrt{17})\left(
x+\frac{-13+5\sqrt{17}}{2}
\right)
\]
of conductor (2).
\end{enumerate}

\medskip

Since $d \neq 5, 17$ in the statement of Theorem~\ref{thm:main} we only have to complete the proof for $d=3$. 
To do this, we must discard the isomorphism $\overline{\rho}_{E,p} \sim \overline{\rho}_{E^\prime,p}$, 
where $E^\prime$ is given in (i) above. The elliptic curve
$E^\prime$ has $j$-invariant $j^\prime=54000$ and so potentially good
reduction at $\mP=(1+\sqrt{3})$; in particular \cite{kraus2},  the order of $\, \overline{\rho}_{E^\prime,p}(I_\mP)$
is $1$,
$2$, $3$, $4$, $6$, $8$ or $24$.
On the other hand, the Frey curve $E$ has potentially  
multiplicative reduction, and $p \nmid \ord_\mP(j)$.
By the theory of the Tate curve \cite[Proposition V.6.1]{SilvermanII} we have $p \mid \overline{\rho}_{E^\prime,p}(I_\mP)$,
giving a contradiction as $p \geq 17$.

\section{Proof of Theorem~\ref{thm:d=17}}

To complete the proof of Theorem~\ref{thm:d=17} it remains to
discard the isomorphism $\overline{\rho}_{E,p} \sim \overline{\rho}_{W,p}$
where $W$ is given in (iii) above. We apply Lemma~1.6 of \cite{HK}---this
is proved for $K=\Q$, but the proof for $K$ a general number field is
identical.
Let $\mP_1$, $\mP_2$ be as in Table~\ref{table:Neven} for $d=17$. 
The curve $E$ has multiplicative reduction at $\mP_i$,
and the valuations of the minimal discriminants are $-8+2p t_i$, where $t_1$, $t_2$
are positive integers.
The curve $W$ has conductor $(2) = \mP_1 \mP_2$ and its minimal
discriminant $\Delta_W$ satisfies $\ord_{\mP_1}(\Delta_W) = 4$ and $\ord_{\mP_2}(\Delta_W) = 2$.
The quantity
\[
 \frac{\ord_{\mP_1}(\Delta_E)\ord_{\mP_2}(\Delta_E)}{\ord_{\mP_1}(\Delta_W)\ord_{\mP_2}(\Delta_W)} = 
\frac{(-8 + 2p t_1)(-8+2p t_2)}{4 \cdot 2} \equiv 8 \pmod{p} 
\]
is a square modulo $p$ if and only if $p \equiv 1$, $7 \pmod{8}$.
It follows from \cite[Lemma 1.6]{HK} that $\overline{\rho}_{E,p} \sim \overline{\rho}_{W,p}$
cannot hold if $p \equiv 3, 5 \pmod{8}$, concluding the proof.

\section{Computational Remarks} \label{sec:end}
In the introduction we indicated that the above strategy can
be applied over other totally real fields (assuming the 
modularity of the Frey curve). However, the computation of
newforms will often be impractical, particularly if 
the levels predicted by level lowering have large norm.
These levels depend crucially on a choice of
odd prime ideal representatives $\cH$ for $\Cl(K)/\Cl(K)^2$. 
In this section we illustrate these computational issues by looking at $K=\Q(\sqrt{30})$
and $K=\Q(\sqrt{79})$.

\bigskip

Let $K=\Q(\sqrt{30})$.
Here $\Cl(K)$ has order $2$, and we can take
$\cH=\{ 1\cdot \OO_K, \, \fp\}$, where $\fp$ is the unique prime
above $3$.
 By computations similar to those leading to
Corollary~\ref{cor:cond}, we obtain four possible levels $\cN_p$. One of
these is
$\cN_p=\mP^8 \cdot \fp^2$,
where $\mP$ is the unique prime above $2$. 
The dimension of the space of cusp forms of level
$\cN_p$ is $26108$, making the computation of newforms 
infeasible with the current {\tt Magma} implementation.

\bigskip

Let $K=\Q(\sqrt{79})$.  Here $\Cl(K)$ has order $3$,
and thus $\Cl(K)/\Cl(K)^2$ is trivial; this is the smallest real
quadratic field for which $\Cl(K)$ and $\Cl(K)/\Cl(K)^2$ differ.
By definition
$\cH=\{ 1 \cdot \OO_K \}$.
We can show by  variants of the arguments in
Section~\ref{sec:irred} that $\overline{\rho}_{E,p}$
is irreducible for $p \ge 17$. Moreover the predicted
levels $\cN_p$ are $\mP$ and $\mP^4$, where $\mP \mid 2$. The dimensions
of the corresponding spaces of cusp forms are
$156$ and $1077$ respectively. Here it feasible to
compute the newforms, and similar arguments to Section~\ref{sec:proof1}
allow us to deduce the following.
\begin{thm}
The Fermat equation~\eqref{eqn:Fermatn} has only trivial solutions
over $K=\Q(\sqrt{79})$, for $n \ge 4$.
\end{thm}

\bigskip

In \cite{Fermat}, the reader will also find recipes for the
possible levels $\cN_p$. The objectives of \cite{Fermat}
are theoretical and there is no need to make the levels $\cN_p$
particularly small. The purpose of the following remarks
is to illustrate the value of Sections~\ref{sec:odd}
and~\ref{sec:scaling} of the current paper, were a finer
analysis of the levels and the effect of scaling the solution
is carried out.
In \cite{Fermat}, the set $\cH$ is taken to be representatives
of $\Cl(K)$ rather than representatives for $\Cl(K)/\Cl(K)^2$.
For the fields $K$ appearing in Theorem~\ref{thm:main}, all
class groups are either trivial or cyclic of order $2$. 
Therefore there is no difference between $\Cl(K)$ and $\Cl(K)/\Cl(K)^2$.
For these fields, the main improvement of the current paper
lies in Section~\ref{sec:scaling} which
radically reduces the possibilities for the even part of the level.
However, to extend the computations to other fields, the
distinction between $\Cl(K)$ and $\Cl(K)/\Cl(K)^2$ becomes crucial.
For example, for $K=\Q(\sqrt{79})$, a set of odd representatives
for $\Cl(K)$ is $\{ 1 \cdot \OO_K, \, \fp_1, \, \fp_2\}$
where $\fp_1 \fp_2=3 \cdot \OO_K$. Following the 
recipe in \cite{Fermat}, the odd part of the level
will be $1 \cdot \OO_K$, $\fp_1^2$ or $\fp_2^2$.
Thus the possibilities for $\cN_p$ include 
$\mP^4 \cdot \fp_1^2$. The dimension of the space
of cusp forms for this level is $12090$, which makes the computation
of newforms impractical. 
Finally, we point out that the even part of the level
given by the recipe in \cite{Fermat} can be as large as $\mP^{10}$. Even if the odd
part of the level is taken to be trivial, the dimension of 
the space of cusp forms of level $\mP^{10}$ is $64596$, which again
is too large. It is clear that the refinements of 
Sections~\ref{sec:odd} and~\ref{sec:scaling} are
required for $K=\Q(\sqrt{79})$, and will be needed
if the computations of the current paper are to be extended
to other totally real fields.

\end{document}